\documentclass[11pt]{article}
\usepackage{amssymb}
\usepackage{amsmath}
\usepackage[dvips]{graphics}
\usepackage{epsfig}

\begin{document}
\title{Non-self-adjoint  Hill operators whose spectrum is a real interval}

\author{
\textbf{Vassilis G. Papanicolaou}
\\\\
Department of Mathematics
\\
National Technical University of Athens
\\
Zografou Campus, 157 80, Athens, GREECE
\\
\underline{\tt papanico@math.ntua.gr}}
\maketitle

\begin{abstract}
Let $H = -d^2/dx^2 + q(x)$, $x \in \mathbb{R}$, where $q(x)$ is a periodic potential, and suppose that the spectrum $\sigma(H)$ of $H$ is the positive semi-axis
$[0, \infty)$. In the case where $q(x)$ is real-valued (and locally square-integrable) a well-known result of G. Borg states that $q(x)$ must vanish almost
everywhere. However, as it was first observed by M.G. Gasymov, there is an abundance of complex-valued potentials for which $\sigma(H) = [0, \infty)$.

In this article we conjecture a characterization of all entire complex-valued potentials whose spectrum is $[0, \infty)$. We also present an analog of Borg's result
for complex potentials.
\end{abstract}


\bigskip

{\bf MSC2020 Mathematical Subject Classification.} 34B30; 34L40; 47E05; 47A10.

\section{The complex Hill operator}
Consider the operator $H$ is acting in $L^2(\mathbb{R})$ defined as
\begin{equation*}
Hy = -y'' +  q(x) y,
\qquad
x \in \mathbb{R},
\end{equation*}
where $q(x)$ is \textbf{complex-valued} and $2\pi$-periodic:
\begin{equation*}
q(x + 2\pi) = q(x),
\qquad
x \in \mathbb{R}.
\end{equation*}
If $q(x)$ is \textbf{real-valued} (and locally square-integrable), then it is well known that $H$ is \textbf{self-adjoint}.

There is a huge amount of literature devoted to the self-adjoint case.

The case of a \textbf{complex-valued} $q(x)$ is mathematically intriguing and has been studied extensively too (see, e.g., [1--22] as well as the references therein).
As expected, the theory is quite different from the self-adjoint case.

The recent emergence of the $\mathcal{PT}$-Symmetric Quantum Theory (see, e.g., \cite{B}) provides another strong motivation for studying non-self-adjoint Schrö-
dinger operators (“non-Hermitian Hamiltonians” in the physicists’ terminology), especially in the case where their spectra are real.

\section{Floquet theory, discriminant and spectrum}
Consider the problem
\begin{equation}
H y = -y'' + q(x) y = \lambda y = k^2 y,
\qquad
x \in \mathbb{R},
\label{Ta1}
\end{equation}
where
\begin{equation*}
\lambda = k^2 \in \mathbb{C}
\end{equation*}
is the spectral parameter.


Let $u(x) = u(x; \lambda)$ and $v(x) = v(v; \lambda)$ be the solutions of \eqref{Ta1} such that
\begin{equation*}
u(0; \lambda) = 1, \quad u'(0; \lambda) = 0,
\qquad\qquad
v(0; \lambda) = 0, \quad v'(0; \lambda) = 1,
\end{equation*}
where primes denote derivatives with respect to $x$.


The Wronskian of $u(x)$ and $v(x)$ is identically equal to $1$. In particular, $u(x)$ and $v(x)$ are linearly independent.


Since we have smooth dependence on the parameter $\lambda$, the solutions $u(x; \lambda)$ and $v(x; \lambda)$ are entire in $\lambda$.
Their orders are $\leq 1/2$ \cite{P-T}.


In the case $q(x) \equiv 0$ (the \textbf{unperturbed} case) we have
\begin{equation*}
\tilde{u}(x; \lambda) = \cos\left(\sqrt{\lambda}\, x\right)
\qquad \text{and} \qquad
\tilde{v}(x; \lambda) = \frac{\sin\left(\sqrt{\lambda}\, x\right)}{\sqrt{\lambda}} 
\end{equation*}
(tilded quantities will be associated with the unpertubed case). 

Now, let $\mathcal{S}$ be the ``shift" or monodromy operator
\begin{equation*}
(\mathcal{S}f)(x) = f(x + 2\pi).
\end{equation*}
The periodicity of $q(x)$ implies that the linear operator $\mathcal{S}$ maps solutions of \eqref{Ta1} to
solutions of \eqref{Ta1} for the same value of $\lambda$ (in other words, $\mathcal{S}$ commutes with $H$), and by exploiting this simple observation
one can develop the Floquet/spectral theory of $H$.


For each $\lambda \in \mathbb{C}$ let $\mathcal{W}  = \mathcal{W}(\lambda)$ be the two-dimensional vector space of the solutions of \eqref{Ta1}. The matrix of the operator $\mathcal{S}|_{\mathcal{W}}$ with respect to the basis $(u, v)$ is
\begin{equation}
S = S(\lambda) =
\left[
  \begin{array}{cc}
    u(2\pi; \lambda) & \ v(2\pi; \lambda) \\
    u'(2\pi; \lambda) & \ v'(2\pi; \lambda) \\
  \end{array}
\right]
\label{B5}
\end{equation}
(the matrix $S$ and the vector space $\mathcal{W}$ depend on $\lambda$).

$S$ is the \textbf{Floquet or monodromy matrix} associated with equation \eqref{Ta1} and
\begin{equation*}
\det S(\lambda) \equiv 1.
\end{equation*}
It follows that the characteristic polynomial of $S(\lambda)$ is
\begin{equation*}
\det \left(S - \rho I\right) = \rho^2 - \Delta(\lambda) \, \rho + 1,
\end{equation*}
where 
\begin{equation*}
\Delta(\lambda) = \text{tr} S(\lambda) = u(2\pi; \lambda) + v'(2\pi; \lambda)
\end{equation*}
is the \textbf{Hill discriminant} (also known as \textbf{Lyapunov's function}) of $H$. Actually, $\Delta(\lambda)$ is entire of order $1/2$ \cite{P-T}.

Sometimes we may find more convenient, instead of $\lambda$, to work with the parameter $k$ (recall that $\lambda = k^2$) and, to avoid confusion, whenever we view
the discriminant as a function of $k$, we will denote it by $D(k)$, so that
\begin{equation*}
D(k) = \Delta(k^2) = \Delta(\lambda).
\end{equation*}
Clearly, $D(k)$ is an even entire function of order $1$.


A remarkable result of V.A. Tkachenko \cite{T3} is that for a function $D(k)$ to be the
Hill discriminant of some Hill operator with a $2\pi$-periodic potential $q(x) \in L^2_{\text{loc}}(\mathbb{R})$, it is \textbf{necessary and
sufficient} that it be an even entire function (of order $1$) of exponential type $2\pi$, which may be represented in the form
\begin{equation}
D(k) = 2\cos(2\pi k) + 2\pi \langle q\rangle \frac{\sin(2\pi k)}{k} - \pi^2 \langle q\rangle^2  \frac{\cos(2\pi k)}{k^2} + \frac{h(k)}{k^2},
\qquad
k \in \mathbb{C},
\label{CC1}
\end{equation}
where
\begin{equation}
\langle q\rangle = \frac{1}{2\pi} \int_0^{2\pi} q(x) \, dx
\label{CC2}
\end{equation}
and $h(k)$ is an (even) entire function of order $\leq 1$; if the order of $h(k)$ is $1$, then its type is $\leq 2\pi$. Furthermore, $h(k)$ satisfies the conditions
\begin{equation}
\int_{-\infty}^{\infty} \left|h(k)\right|^2 dk < \infty
\qquad \text{and} \qquad
\sum_{n = -\infty}^{\infty} \left|h\left(\frac{n}{2}\right)\right| < \infty.
\label{CC3}
\end{equation}
Incidentally, let us mention that in the discrete case, where the corresponding discrete complex Hill operator is acting on $\ell^2(\mathbb{Z})$, $\mathbb{Z}$ being
the integer lattice, the Hill discriminant $\Delta(\lambda)$ associated with an $N$-periodic (discrete) potential (where $N \geq 1$ is an integer) can be any polynomial 
whose leading term is $(-1)^N \lambda^N$. Furthermore, any such polynomial is the discriminant of at least $1$ and at most $N!$ discrete Hill operators \cite{Pa2}.

\medskip

Now, let $\rho_1(\lambda)$ and $\rho_2(\lambda) = \rho_1(\lambda)^{-1}$ be the eigenvalues of $S(\lambda)$, namely the \textbf{Floquet multipliers} of $H$.
We have
\begin{equation*}
\rho_1(\lambda) + \rho_2(\lambda) = \text{tr} S(\lambda) = \Delta(\lambda),
\end{equation*}
and
\begin{equation*}
\rho_1(\lambda), \rho_2(\lambda) = \frac{\Delta(\lambda) \pm\sqrt{\Delta(\lambda)^2 - 4}}{2}.
\end{equation*}
The eigenvectors of $S(\lambda)$ associated with its eigenvalues $\rho_1(\lambda)$ and $\rho_2(\lambda)$ correspond to the the \textbf{Floquet solutions}
$\phi_1(x)$ and $\phi_2(x)$ of \eqref{Ta1} satisfying
\begin{equation*}
\phi_j(x + 2\pi) = (\mathcal{S} \phi_j)(x) = \rho_j \phi_j(x),
\qquad
j= 1,2.
\end{equation*}


Notice that $\rho_1(\lambda) = \rho_2(\lambda)$ can happen only if $\rho_1(\lambda) = \rho_2(\lambda) = \pm1$ (equivalently, $\Delta(\lambda) = \pm2$).
In this case we may not have two linearly independent Floquet solutions. If it happens that for a given $\lambda$ satisfying $\Delta(\lambda) = \pm2$ 
two linearly independent Floquet solutions exist, then we say we have \textbf{coexistence}.

It is sometimes more convenient to view $\rho_1(\lambda)$ and $\rho_2(\lambda)$ as the two branches of a (single-valued) analytic function $\rho(\lambda)$
defined on the Riemann surface of the function $\sqrt{\Delta(\lambda)^2 - 4}$ (generically, this Riemann surface is not compact since $\Delta(\lambda)^2 - 4$
is entire of order $1/2$ and, consequently, it has infinitely many zeros by the Hadamard Factorization Theorem \cite{A}; unless all but finitely many zeros of
$\Delta(\lambda)^2 - 4$ have even multiplicity, the Riemann surface is not compact).
Thus,
\begin{equation*}
\rho(\lambda) + \frac{1}{\rho(\lambda)} = \Delta(\lambda),
\qquad
\rho(\lambda) = \frac{\Delta(\lambda) + \sqrt{\Delta(\lambda)^2 - 4}}{2}
\end{equation*}
and $\rho(\lambda)$ can be called the \textbf{Floquet multiplier} associated with \eqref{Ta1}.

The fact that $\Delta(\lambda)$ is entire implies that $\rho(\lambda)$ has neither zeros nor poles (nor essential singularities) for any finite $\lambda$.
Therefore, the only possible singularities of $\rho(\lambda)$ are square-root branch points at which we must necessarily have
$\rho(\lambda) = \pm 1$ (equivalently, $\Delta(\lambda) = \pm2$).

Actually, $\rho(\lambda)$ must have at least one branch point, since if it had no branch points,
then it would have been an entire function of order $\leq 1/2$ with no zeros, therefore a constant (by the Hadamard Factorization Theorem \cite{A}), which is impossible since $\Delta(\lambda)$ is not a constant.

In some sense, $\rho(\lambda)$ can be viewed as the analog of the exponential function for the Riemann surface of $\sqrt{\Delta(\lambda)^2 - 4}$. Also,
\begin{equation*}
\left[\log\rho(\lambda)\right]' = \frac{\rho'(\lambda)}{\rho(\lambda)} = \frac{\Delta'(\lambda)}{\sqrt{\Delta(\lambda)^2 - 4}}
\end{equation*}
and, since $\rho(\lambda)$ is single-valued on the Riemann surface, we have that the holomorphic differential
\begin{equation*}
\frac{\Delta'(\lambda)}{\sqrt{\Delta(\lambda)^2 - 4}} \, d\lambda
\end{equation*}
has period $2\pi i$ ($\log\rho(\lambda)$ is the \textbf{Floquet exponent}).

The values of $\lambda$ for which $\rho(\lambda) = 1$ (equivalently, $\Delta(\lambda) = 2$) are the \textbf{periodic eigenvalues}
of $H$, since, in this case, any associated Floquet solution is $2\pi$-periodic.


The values of $\lambda$ for which $\rho(\lambda) = -1$ (equivalently, $\Delta(\lambda) = -2$) are the \textbf{antiperiodic eigenvalues} of $H$, since, in this case,
any associated Floquet solution is $2\pi$-\textbf{antiperiodic}, namely
\begin{equation*}
\phi(x + 2\pi) = -\phi(x).
\end{equation*}

As we have already mentioned, $S(\lambda)$ can have a Jordan anomaly only if $\rho(\lambda) = \pm 1$ (equivalently, only if $\Delta(\lambda) = \pm 2$) and in the presence of such an anomaly the matrix $S(\lambda)$ is similar to the Jordan canonical matrix
\begin{equation*}
\left[
  \begin{array}{cc}
    \pm 1 & 1 \\
    0 & \pm 1 \\
  \end{array}
\right].
\end{equation*}
Let us mention that $\lambda^{\star}$ can be a zero of $\Delta(\lambda)^2 - 4$ of even multiplicity, so that $\lambda^{\star}$ is not a branch point of $\rho(\lambda)$,
and, yet, $S(\lambda^{\star})$ may not be diagonalizable. If this is the case, we say that the Floquet matrix $S(\lambda)$ has a \emph{pathology of the second
kind at} $\lambda^{\star}$.

If for some $\lambda = \lambda^{\star}$ we have coexistence of two periodic or, respectively, antiperiodic solutions, then
\begin{equation*}
S(\lambda^{\star}) = \left[
  \begin{array}{cc}
     1 & 0 \\
    0 &  1 \\
  \end{array}
\right],
\quad
\text{respectively}\quad
S(\lambda^{\star}) = \left[
  \begin{array}{cc}
    -1 & 0 \\
    0 & - 1 \\
  \end{array}
\right].
\end{equation*}

If $\lambda^{\star}$ is a periodic eigenvalue for which we have coexistence of two periodic solutions, then $\lambda^{\star}$ is a zero of $\Delta(\lambda) - 2$
of multiplicity $\geq 2$. Likewise, if $\lambda^{\star}$ is an antiperiodic eigenvalue for which we have coexistence of two antiperiodic solutions, then 
$\lambda^{\star}$ is a zero of $\Delta(\lambda) + 2$ of multiplicity $\geq 2$ (in this sense we may say that the algebraic multiplicity of a periodic/antiperiodic
eigenvalue is greater or equal to its geometric multiplicity).


The last statement follows from the formula (which can be derived by writing \eqref{Ta1} for $u(x; \lambda)$ and $v(x; \lambda)$, then differentiating with respect to
$\lambda$ and applying variation of parameters)
\begin{align}
\Delta'(\lambda) =&\, u(2\pi; \lambda) \int_0^{2\pi} u(x; \lambda) v(x; \lambda) dx
-v(2\pi; \lambda) \int_0^{2\pi} u(x; \lambda)^2 dx
\nonumber
\\
& +u'(2\pi; \lambda) \int_0^{2\pi} v(x; \lambda)^2 dx
-v'(2\pi; \lambda) \int_0^{2\pi} u(x; \lambda) v(x; \lambda) dx.
\nonumber
\end{align}

\subsection{The spectrum}
The spectrum $\sigma(H)$ of $H$ is characterized as
\begin{align}
\sigma(H) &= \{\lambda \in \mathbb{C} : |\rho(\lambda)| = 1\}
= \{\lambda \in \mathbb{C} : \rho(\lambda) = e^{i \theta}, \quad 0 \leq \theta \leq \pi\}
\nonumber
\\
 &= \{\lambda \in \mathbb{C} : \Delta(\lambda) \in [-2, 2]\}
= \{\lambda \in \mathbb{C} : \Delta(\lambda) = 2\cos \theta, \quad 0 \leq \theta \leq \pi\}.
\nonumber
\end{align}
Notice that $\sigma(H)$ is an unbounded closed subset of $\mathbb{C}$ (this follows, e.g., from the fact that $\Delta(\lambda)$ is entire of order
$1/2$ and, consequently, by the Hadamard Factorization Theorem \cite{A} it takes every value in $[-2, 2]$ infinitely many times).


More precisely [5, 12, 14, 16, 19],
$\sigma(H)$ is a countable
system (i.e. union) of analytic arcs, where the analyticity of such an arc may fail only at a point $\lambda$ such that $\Delta'(\lambda) = 0$ (while
$\Delta(\lambda) = 2 \cos\theta$ for some $\theta \in [0, \pi]$, so that $\lambda$ lies in the spectrum). Furthermore, the
resolvent set $\mathbb{C} \smallsetminus \sigma(H)$ of $H$ is path-connected. In particular, $\sigma(H)$ cannot contain closed curves and, also, it
cannot be a piecewise analytic curve without an endpoint. Asymptotically, the spectral arcs approach the half-line (the asymptotic form of the spectrum)
\begin{equation*}
\ell_{\langle q\rangle} = \{z \in \mathbb{C} \,:\, z = \langle q\rangle + x, \ x\geq 0\}.
\end{equation*}

A rather trivial observation is that if $\lambda^{\star}$ is a periodic or antiperiodic eigenvalue, then  $\Delta(\lambda^{\star}) = \pm2$,
hence $\lambda^{\star} \in \sigma(H)$.

\section{The case where $\sigma(H)$ is a single analytic arc}
Suppose that the spectrum $\sigma(H)$ is an analytic (connected) curve. Since $\mathbb{C} \smallsetminus \sigma(H)$ is path-connected, $\sigma(H)$ must have one
(and only one) endpoint, say $\lambda_0$.

By replacing $q(x)$ by $q(x) - \lambda_0$, we can assume that the endpoint of $\sigma(H)$ is $0$. 

Suppose $\Delta(\lambda^{\star})^2 - 4 = 0$. Then $\lambda^{\star} \in \sigma(H)$. Let us assume that $\lambda^{\star} \neq 0$ so that $\lambda^{\star}$ is an
``interior" point of $\sigma(H)$. From the Taylor expansion of $\Delta(\lambda)$ about $\lambda^{\star}$ we get
\begin{equation*}
\Delta(\lambda) = \pm2 + c (\lambda - \lambda^{\star})^d + O\left[(\lambda - \lambda^{\star})^{d+1}\right],
\qquad
\lambda \to \lambda^{\star},
\end{equation*}
where $d$ is an integer $\geq 1$ and $c \neq 0$. Then, the assumption that $\lambda^{\star}$ is an interior point of $\sigma(H)$ forces $d = 2$ Hence, $\lambda^{\star}$
cannot be a branch point of $\rho(\lambda)$.

It follows that $0$ is the unique branch point of $\rho(\lambda)$. Thus,
\begin{equation*}
\rho(\lambda) = f\left(\sqrt{\lambda}\right) = f(k)
\qquad (\text{since }\;\lambda = k^2),
\end{equation*}
where $f(k)$ is entire of order $1$ and has no zeros. Furthermore $0$ is a branch point of $\rho(\lambda)$ and, hence, $\rho(0) = \pm 1$.
Therefore, $\rho(\lambda)$ must be of the form
\begin{equation*}
\rho(\lambda) = \pm e^{i \alpha \sqrt{\lambda}},
\end{equation*}
where $\alpha \ne 0$ is a complex constant.

Hence,
\begin{equation*}
\Delta(\lambda) = \rho(\lambda) + \rho(\lambda)^{-1} = \pm 2\cos\left(\alpha \sqrt{\lambda}\right),
\end{equation*}
and the general characterization of the discriminant given in \eqref{CC1} implies that $\alpha = 2\pi$ and $\langle q - \lambda_0\rangle = 0$ (i.e. for our original 
$q(x)$ we must have $\langle q\rangle = \lambda_0$).

Furthermore, again by \eqref{CC1}, we must have
\begin{equation*}
\Delta(\lambda) = 2\cos\left(2\pi\sqrt{\lambda}\right),
\qquad \text{hence} \quad
\rho(\lambda) = e^{2\pi i\sqrt{\lambda}}
\end{equation*}
and, consequently,
\begin{equation*}
\sigma(H) = [0, \infty)
\end{equation*}
(for our original $q(x)$ we must have $\sigma(H) = \langle q\rangle + [0, \infty)$).
Notice also that $\rho(0) = 1$, hence $0$ is a periodic eigenvalue. Furthermore, $\Delta'(\lambda) = -2\pi \sin\left(2\pi\sqrt{\lambda}\right) / \sqrt{\lambda}$,
hence $\Delta'(0) = -4\pi^2 \neq 0$, which implies that for $\lambda = 0$ we cannot have coexistence.


Thus, $S(\lambda)$ does not have a \emph{pathology of the first kind}
at $\lambda = 0$ (a pathology of the first kind at $\lambda^{\star}$ occurs if $\lambda^{\star}$ is a 
branch point of $\rho(\lambda)$ and at the same time we have coexistence of two periodic or antiperiodic solutions at $\lambda = \lambda^{\star}$).

\section{The self-adjoint case}
In the \textbf{self-adjoint case} (i.e. when $q(x)$ is real-valued) $\lambda^{\star}$ is a double zero of $\Delta(\lambda) - 2$ if and only if we have coexistence
of periodic solutions for $\lambda = \lambda^{\star}$, while $\lambda^{\star}$ is a double zero of $\Delta(\lambda) + 2$ if and only if we have coexistence
of antiperiodic solutions for $\lambda = \lambda^{\star}$. Furthermore, $\Delta(\lambda)^2 - 4$ does not have any zeros with multiplicity $> 2$.  In this sense, 
algebraic multiplicity equals geometric multiplicity. Also, a point $\lambda^{\star}$ is a branch point of the Floquet multiplier $\rho(\lambda)$ if and only if $S(\lambda^{\star})$ has a Jordan anomaly.

The spectrum is a union of closed intervals (the \textbf{bands}) separated by open intervals (the \textbf{gaps}):
\begin{equation*}
\sigma(H) = \bigcup_{n \geq 0} \left[\lambda_{2n}, \lambda_{2n+1}\right],
\qquad
\lambda_0 < \lambda_1 \leq \lambda_2 < \lambda_3 \leq \lambda_4 < \lambda_5 \leq \lambda_6 < \cdots
\end{equation*}

$\lambda_0$ and $\lambda_{4j-1} \leq \lambda_{4j}$, $j \leq 1$ are the periodic eigenvalues, while
$\lambda_{4j-3} \leq \lambda_{4j-2}$, $j \leq 1$ are the antiperiodic eigenvalues.

If for some $n \geq 1$ we have that $\lambda_{2n-1} = \lambda_{2n}$,
then the corresponding gap
$(\lambda_{2n-1}, \lambda_{2n})$ of the spectrum is \textbf{closed} (i.e. empty) and we have coexistence of two 
linearly independent periodic or antiperiodic solutions.

If $\lambda_{2n-1} < \lambda_{2n}$, then there is no coexistence neither at $\lambda_{2n-1}$ nor at $\lambda_{2n}$.

Clearly, the Dirichlet spectrum $\{\mu_1, \mu_2, \ldots\}$ of $H$ on the interval $(0, 2\pi)$ coincides with the set of (distinct) zeros of the entire function
$v(2\pi; \lambda)$.

In the self-adjoint case all the zeros of $v(2\pi; \lambda)$ are simple and, of course, real. Furthermore, if $v(2\pi; \mu) = 0$, then the Floquet matrix at
$\lambda = \mu$ becomes
\begin{equation*}
S(\mu) =
\left[
  \begin{array}{cc}
    u(2\pi; \mu) & \ 0 \\
    u'(2\pi; \mu) & \ v'(2\pi; \mu) \\
  \end{array}
\right],
\end{equation*}
hence the \textbf{real} quantities $u(2\pi; \mu)$ and $u'(2\pi; \mu)$ are the eigenvalues of $S(\mu)$, i.e. the Floquet multipliers.
In particular, $u(2\pi; \mu) u'(2\pi; \mu) = 1$ and, consequently
\begin{equation*}
|\Delta(\mu)| = |u(2\pi; \mu) + u'(2\pi; \mu)| = |u(2\pi; \mu)| + |u'(2\pi; \mu)| \geq 2. 
\end{equation*}

Actually, we have
\begin{equation*}
\lambda_0 < \lambda_1 \leq \mu_1 \leq \lambda_2 < \lambda_3 \leq \mu_2 \leq \lambda_4 < \lambda_5 \leq \mu_3 \leq \lambda_6 < \cdots.
\end{equation*}


There is a very short proof of all the above properties of the self-adjoint case. First we check them for the trivial case $q(x) \equiv 0$ and then we
consider the continuous deformation of potentials
\begin{equation*}
t q(x),
\qquad
0 \leq t \leq 1,
\end{equation*}
and exploit the continuous dependence on $t$ (notice that, by self-adjointness all motion of the $\lambda$'s and $\mu$'s is confined on the real axis).

\section{A well-known theorem of G. Borg}

In his famous paper \cite{Bo} (see also \cite{U}) among many other inverse spectral results regarding the Sturm-Liouville operator, Borg has shown that for a \textbf{real-valued} potential $q(x) \in L^2_{loc}(\mathbb{R})$:

If  $\sigma(H) = [0, \infty)$,   then   $q(x) = 0$ a.e.

Actually, Borg proved a more general statement. He showed that if all the gaps corresponding to antiperiodic eigenvalues are closed, then
\begin{equation*}
q(x + \pi) = q(x)\ \text{a.e.}
\end{equation*}

\medskip

\underline{\textbf{QUESTION}}\textbf{:} Are there analogs or extensions to Borg’s theorem in the complex potential case?

\medskip

It is worth mentioning that Borg's theorem fails in the case where the potential $q(x)$ is quasi-periodic (we believe that Borg's theorem also fails in the case
where $q(x)$ is limit-periodic).


\section{M.G. Gasymov's discovery}
The case of a nonreal $q(x)$, however, is quite different. Gasymov \cite{G} made the remarkable discovery that if
\begin{equation*}
q(x) = \sum_{n=1}^{\infty} B_n e^{inx},
\qquad \text{with }\;
\sum_{n=1}^{\infty} |B_n| < \infty,
\end{equation*}
then the equation
\begin{equation*}
H y = -y'' + q(x) y = k^2 y,
\end{equation*}
has a Floquet solution of the form
\begin{equation*}
\phi(x; k) = e^{ikx}\left(1 + \sum_{n=1}^{\infty} \frac{1}{n + 2k} \sum_{\ell = n}^{\infty} c_{n\ell} e^{i\ell x}\right),
\end{equation*}
where \textbf{the coefficients} $c_{n\ell}$ \textbf{do not depend on} $k$ and satisfy
\begin{equation*}
\sum_{n=1}^{\infty} \frac{1}{n} \sum_{\ell = n+1}^{\infty} \ell (\ell - n) |c_{n\ell}| < \infty
\qquad\text{and}\qquad
\sum_{n=1}^{\infty} n |c_{n\ell}| < \infty.
\end{equation*}
It follows that the Floquet multiplier is
\begin{equation*}
e^{2\pi i k} = e^{2\pi i\sqrt{\lambda}},
\end{equation*}
and consequently, $\sigma(H) = [\,0, \infty)$.

Notice also that $\phi(x; k)$ is meromorphic in $k$ and its poles are simple. Furthermore, every pole is of the form $-n/2$, where
$n$ is a positive integer, and if $k \neq -n/2$, $n = 0, 1, \ldots$, then $\phi(x; -k)$ is the other Floquet solution.

Actually, since the spectral properties of the operator $H$ depend continuously on $q(x)$ with respect to the $L^2(0, 2\pi)$-norm, it follows that
for the weaker assumption $\sum_{n=1}^{\infty} |B_n|^2 < \infty$ we still have
$\sigma(H) = [\,0, \infty)$ \cite{Sh}.

It is also worth mentioning that there are multidimensional analogs of Gasymov's result (see, e.g., \cite{K}).

Since
\begin{equation*}
\Delta(\lambda) = 2\cos(2\pi \sqrt{\lambda})
\qquad
\Rightarrow
\qquad
\Delta(\lambda)^2 - 4 = -4 \sin^2(2\pi \sqrt{\lambda}),
\end{equation*}
the zeros of $\Delta(\lambda)^2 - 4$ are (counting multiplicities)
\begin{equation*}
\left(\frac{n}{2}\right)^2,
\qquad
n \in \mathbb{Z}.
\end{equation*}
Notice that $0$ is a simple zero of $\Delta(\lambda)^2 - 4$, while all other zeros, namely the zeros $n^2/4$, $n \geq 1$, are
double.

Clearly, the only branch point of the Floquet multiplier $\rho(\lambda) =  e^{2\pi i\sqrt{\lambda}}$ is
$\lambda = 0$. However, $S(n^2/4)$, where $S(\lambda)$ is given by \eqref{B5}, may not be diagonalizable for nonzero values of $n$ (pathology of the second kind).

There is an easy way to (partly) understand Gasymov's result. In the equation
\begin{equation*}
-y'' + q(x) y = k^2 y,
\qquad\qquad
q(x) = \sum_{n=1}^{\infty} B_n e^{inx},
\end{equation*}
we substitute
\begin{equation*}
z = e^{ix},
\qquad
w(z) = w(e^{ix}) = y(x).
\end{equation*}
Then, the equation becomes
\begin{equation*}
z^2 w''(z) + z w'(z) + P(z) w(z) = k^2 w(z),
\qquad\text{with}\quad
P(z) = \sum_{n=1}^{\infty} B_n z^n.
\end{equation*}
This equation has a regular singular point at $z = 0$. Therefore its solutions can be expressed in Frobenius series.
The indicial equation is
\begin{equation*}
r^2 = k^2,
\qquad\text{thus}\quad
r = \pm k,
\end{equation*}
and, hence, the Frobenius solutions are (at least for $k \neq n/2$, $n = 0, \pm1, \ldots$)
\begin{equation*}
w(z) = z^{\pm k}\sum_{n = 0}^{\infty} a_n z^n,
\end{equation*}
which implies that the Floquet multiplier of the original equation is $e^{2\pi i k}$ and, consequently, the spectrum is $[0, \infty)$.

\section{An example}
For a fixed integer $m \geq 1$ and a fixed complex number $a \neq 0$, with $|a| \neq 1$, we set
\begin{equation}
q_m(x) = \frac{2 m^2 a e^{imx}}{(a e^{imx} + 1)^2} = \frac{2 m^2 a^{-1} e^{-imx}}{(a^{-1} e^{-imx} + 1)^2}
= \frac{m^2}{2} \, \text{sech}^2\left(\frac{\xi + imx}{2}\right),
\label{B14}
\end{equation}
where $\xi = \log a$.
Notice that for $|a| < 1$ we have
\begin{equation*}
q_m(x) = 2 m^2 a \sum_{n=1}^{\infty} (-1)^{n+1} n e^{inmx},
\end{equation*}
while for $|a| > 1$ we have
\begin{equation*}
q_m(x) = 2 m^2 a^{-1} \sum_{n=1}^{\infty} (-1)^{n+1} n e^{-inmx}.
\end{equation*}
Then, one Floquet solution of the equation
\begin{equation}
-y'' + q_m(x) y = \lambda y = k^2 y
\label{B15}
\end{equation}
is
\begin{equation}
\phi(x; k) = e^{ikx}\left[1 -  \frac{1}{k + (m/2)} \cdot \frac{m a e^{imx}}{a e^{imx} + 1}\right]
\label{A7cc}
\end{equation}
(in the case $|a| < 1$ this is the Gasymov solution).

Now, unless $k = m/2$, we have that $\phi(x; -k)$ is also a Floquet solution and, furthermore, $\phi(x; k)$ and $\phi(x; -k)$ are linearly independent
for $k \neq 0$ (and $k \neq \pm m/2$). Actually, for $k \neq m/2$ the Wronskian of $\phi(x; k)$ and $\phi(x; -k)$ is $-2ik$. Thus, we have coexistence
for all $k \neq 0, \pm m/2$.

For $k = 0$, i.e. for $\lambda = 0$, another solution is
\begin{equation*}
\left(x -  \frac{4}{im} \cdot \frac{1}{a e^{imx} - 1}\right) \phi(x; 0),
\end{equation*}
which is, obviously, not periodic. Hence, we do not have coexistence. Furthermore, let us notice that $\lambda = 0$ is a simple zero of
\begin{equation*}
\Delta(\lambda)^2 - 4 = -4 \sin^2\left(2\pi \sqrt{\lambda}\right).
\end{equation*}

For $k = \pm m/2$, i.e. for $\lambda = m^2/4$, another solution is
\begin{equation*}
\left(2iamx + a^2 e^{imx} - e^{-imx}\right) \phi(x; m/2),
\end{equation*}
which is, obviously, neither periodic nor antiperiodic. Hence, again, we do not have coexistence. However, $\lambda = m^2/4$ is a
\textbf{double} zero of $\Delta(\lambda)^2 - 4 = -4 \sin^2\left(2\pi \sqrt{\lambda}\right)$ (pathology of the second kind).


The solution $v(x; \lambda)$ of \eqref{B15} satisfying $v(0; \lambda) = 0$ and $v'(0; \lambda) = 1$ is
\begin{equation}
v(x; \lambda) = \frac{1}{8ik (k^2 - m^2/4)} \left[C_m(x\,; k\,; a) e^{ikx} - C_m(x\, ;-k\,; a) e^{-ikx}\right],
\label{B17}
\end{equation}
where
\begin{equation}
C_m(x\, ;k\,; a) = \left(\frac{a-1}{a+1} \, m + 2k\right) \left(m + 2k - \frac{2ma e^{imx}}{a e^{imx} + 1} \right)
\label{B17tt}
\end{equation}
(as usual, $\lambda = k^2$). Formula \eqref{B17} is valid for every $\lambda \in \mathbb{C}$. For instance,
for $\lambda = 0$ formula \eqref{B17} becomes
\begin{equation*}
v(x; 0) = \frac{4a (e^{imx} - 1) + im (a-1) (ae^{imx} - 1) x}{im (a+1) (a e^{imx} + 1)}.
\end{equation*}
For $x = 2\pi$ formula \eqref{B17} yields
\begin{equation}
v(2\pi; \lambda) = \frac{\sin\left(2\pi\sqrt{\lambda}\right)}{\sqrt{\lambda} \, (\lambda - m^2/4)} \left[\lambda - \left(\frac{a-1}{a+1}\right)^2 \frac{m^2}{4} \right].
\label{B17b}
\end{equation}
From formula \eqref{B17b} we see that the zeros of $v(2\pi; \lambda)$ (counting multiplicities) are
\begin{equation*}
\mu_n = \frac{n^2}{4},
\quad
n \geq 1,\ n \neq m;
\qquad \qquad
\mu_m = \left(\frac{a-1}{a+1}\right)^2 \frac{m^2}{4}
\end{equation*}
and, hence, for each $n \geq 1$, $n \neq m$, there are nonzero values of $a$ for which the number $n^2/4$ becomes a double zero of $v(2\pi; \lambda)$.


The solution $u(x; \lambda)$ of \eqref{B15} satisfying $u(0; \lambda) = 1$ and $u'(0; \lambda) = 0$ is
\begin{align}
u(x; \lambda) = &\frac{2(a+1)^2 k^2 + (a^2 -1) mk - 2am^2}{2(a+1)^2 (2k-m) k}\phi(x; k)
\nonumber
\\
+ &\frac{2(a+1)^2 k^2 - (a^2 -1) mk - 2am^2}{2(a+1)^2 (2k+m) k}\phi(x; -k),
\label{B17c}
\end{align}
where $\phi(x; k)$ is given by \eqref{A7cc} (as usual, $\lambda = k^2$). Formula \eqref{B17c} is valid for every $\lambda \in \mathbb{C}$. For instance,
for $\lambda = 0$ formula \eqref{B17c} becomes
\begin{equation*}
u(x; 0) = \frac{a (a^2 + 4a - 1) e^{imx} - 2i a^2 m x e^{imx} + 2i a m x - a^2 + 4a + 1}{(a+1)^2 (a e^{imx} + 1)}.
\end{equation*}
For $x = 2\pi$ formula \eqref{B17c} yields
\begin{equation}
u'(2\pi; \lambda)
= -\frac{\sin\left(2\pi\sqrt{\lambda}\right)}{\sqrt{\lambda} \, (\lambda - m^2/4)}
\left[\lambda^2 - \frac{(a^2 + 6a + 1)m^2}{4(a+1)^2} \, \lambda + \frac{a^2 m^4}{(a+1)^4}\right].
\label{B17bc}
\end{equation}
From formula \eqref{B17bc} we see that the zeros of $u'(2\pi; \lambda)$ (counting multiplicities) are
\begin{equation*}
\nu_n = \frac{n^2}{4},
\
n \geq 1,\ n \neq m;
\
\nu_0,\nu_m = \frac{a^2 + 6a + 1 \pm(a-1) \sqrt{a^2 + 14a + 1}}{4(a+1)^2}\cdot \frac{m^2}{4}
\end{equation*}
(thus, for $a \neq 0, 1$ we get that $\nu_0,\nu_m \neq 0$ and $\nu_0,\nu_m \neq m^2/4$).

As it is well known, the potential $q_m(x)$ of \eqref{B14} is obtained by applying a Darboux
transformation to the trivial potential $q(x) \equiv 0$.

\section{A conjecture}
\textbf{Conjecture.} Let $q(x)$ be an entire and $2\pi$-periodic function of $x$. If the spectrum of the operator
$H = -d^2/dx^2 + q(x)$ is
\begin{equation*}
\sigma(H) = [0, \infty),
\end{equation*}
then
\begin{equation*}
q(x) = \sum_{n=1}^{\infty} A_n e^{-inx}
\qquad \text{or } \qquad
q(x) = \sum_{n=1}^{\infty} B_n e^{inx}.
\end{equation*}

\medskip

\textbf{Terminology.} We call \emph{Gasymov potential} any (not necessarily entire) periodic function $G(x)$ whose Fourier series
expansion contains only positive or only negative frequencies.

\medskip

A small indication in favor of the conjecture is the following:


If the Fourier expansion of $q(x)$ contains both positive and negative frequencies, then the resulting 
equation with respect to $z = e^{ix}$ has a singular singular point at $z = 0$.

\section{The shifted operator}
Let $\xi$ be a given real number. We introduce the shifted operator
\begin{equation*}
(H_{\xi} \,y)(x) = -y''(x) + q_{\xi}(x) \, y(x)
\qquad
\text{acting in }\; L^2(\mathbb{R}),
\end{equation*}
where
\begin{equation*}
q_{\xi}(x) = q(x + \xi)
\end{equation*}
(thus $H_0 = H$).

\medskip

\textbf{Notation.} If $A$ is a quantity associated with the operator $H$, the corresponding quantity associated with the operator $H_{\xi}$ will be
denoted by $A_{\xi}$.

\medskip

Suppose that $\phi(x)$ is a Floquet solution of $Hy = \lambda y$ associated with the Floquet multiplier $\rho(\lambda)$, so that
\begin{equation*}
\phi(x + 2\pi) = \rho(\lambda)  \phi(x).
\end{equation*}
Then, $\phi(x + \xi)$ (as a function of $x$) satisfies the equation $H_{\xi} \, y = \lambda y$ and we also have that
$\phi(x + 2\pi + \xi) = \rho(\lambda) \phi(x + \xi)$, which means that $\phi(x + \xi)$ (as a function of $x$) is a Floquet solution of
$H_{\xi} \, y = \lambda y$ associated with the Floquet multiplier $\rho(\lambda)$. Furthermore, since this is true for every $\lambda \in \mathbb{C}$ it follows
that
\begin{equation}
\rho_{\xi}(\lambda) \equiv \rho(\lambda)
\label{B22}
\end{equation}
i.e. the operators $H$ and $H_{\xi}$ have the same Floquet multiplier and, consequently,
\begin{equation}
\sigma(H_{\xi}) = \sigma(H),
\label{B23}
\end{equation}
thus the spectrum of $H$ remains invariant under the shift by $\xi$.

We also get that
\begin{equation}
\Delta_{\xi}(\lambda) \equiv \Delta(\lambda),
\qquad\text{i.e. }\; u_{\xi}(2\pi; \lambda) + v_{\xi}'(2\pi; \lambda) \equiv u(2\pi; \lambda) + v'(2\pi; \lambda).
\label{B24}
\end{equation}

Suppose now that $q(x)$ is analytic in a strip $\mathcal{T}$ of the form $a < \Im(x) < b$ containing the real axis.
Then $q_{\xi}(x) = q(x + \xi)$
makes sense for $\xi \in \mathcal{T}$ and $x \in \mathbb{R}$. Therefore, by analytic continuation the equations \eqref{B22}, \eqref{B23}, and \eqref{B24}
remain true for all $\xi \in \mathcal{T}$, $x \in \mathbb{R}$. If, in particular, $q(x)$ is entire in $x$, then they remain true for all $\xi \in \mathbb{C}$.

If, however, $q(x)$ is \textbf{meromorphic} in $x$, the equations \eqref{B22}, \eqref{B23}, and \eqref{B24} may not quite hold for every $\xi \in \mathbb{C}$.
For instance, let
\begin{equation*}
q(x) = \frac{e^{ix}}{1 - (1/2) \, e^{ix}}.
\end{equation*}
Clearly, $q(x)$ is meromorphic and
\begin{equation*}
q(x) = \sum_{n=1}^{\infty} \frac{e^{inx}}{2^{n-1}},
\qquad
x \in \mathbb{R}.
\end{equation*}
Thus, $q(x)$ is a Gasymov potential and, consequently, $\sigma(H) = [0, \infty)$. Now, let us consider the shifted potential
\begin{equation*}
q_{\xi}(x) = \frac{e^{i\xi} e^{ix}}{1 - (1/2) \, e^{i\xi} e^{ix}}.
\end{equation*}
By choosing $\xi = -i \log 4$ we get
\begin{equation*}
q_{\xi}(x) = \frac{4 e^{ix}}{1 - 2 e^{ix}} = \frac{-2}{1 - (1/2) \, e^{-ix}}
= -2 -\sum_{n=1}^{\infty} \frac{e^{-inx}}{2^{n-1}},
\qquad
x \in \mathbb{R},
\end{equation*}
from which we see that $q_{\xi}(x) + 2$ is a Gasymov potential and hence
\begin{equation*}
\sigma(H_{\xi}) = [-2, \infty) \neq \sigma(H).
\end{equation*}

\section{Asymptotic formulas}
Suppose $q(x)$ is in $C^2$. Then (see, e.g., \cite{P-T}),
\begin{align}
v(x; \lambda) = &\, \tilde{v}(x; \lambda) - \frac{\cos\left(\sqrt{\lambda}\, x\right)}{2 \sqrt{\lambda}} Q(x)
+ \frac{\tilde{v}(x; \lambda)}{4 \lambda} \left[q(x) + q(0) - \frac{Q(x)^2}{2} \right]
\nonumber
\\
& + O\left(\frac{e^{\left|\Im\left(\sqrt{\lambda}\right)\right| x}}{\left|\lambda\right|^2}\right),
\quad
\lambda \to \infty,
\label{T28a}
\end{align}
where
\begin{equation}
\tilde{v}(x; \lambda) = \frac{\sin\left(\sqrt{\lambda}\, x\right)}{\sqrt{\lambda}} 
\qquad\text{and}\qquad
Q(x) = \int_0^x q(\xi)d\xi
\label{T1}
\end{equation}
(recall that $\tilde{v}(x; \lambda)$ is the corresponding solution of the unperturbed problem).

Thus, if
\begin{equation*}
\langle q\rangle = \frac{Q(2\pi)}{2\pi} = \frac{1}{2\pi}\int_0^{2\pi} q(\xi)d\xi = 0,
\end{equation*}
then \eqref{T28a} implies
\begin{equation}
v(2\pi; \lambda) = \tilde{v}(2\pi; \lambda)
+ \frac{\tilde{v}(2\pi; \lambda)}{2 \lambda}\, q(0)
+ O\left(\frac{e^{2\pi\left|\Im\left(\sqrt{\lambda}\right)\right|}}{\left|\lambda\right|^2}\right),
\qquad
\lambda \to \infty.
\label{T28c}
\end{equation}
If $N$ is a sufficiently large integer, then $v(2\pi; \lambda)$ has exactly $N$ zeros (counting multiplicities) in the open half-plane \cite{P-T}
\begin{equation}
\Re(\lambda) < \left(\frac{N}{2} + \frac{1}{4}\right)^2
\label{T28d}
\end{equation}
(notice that $\tilde{v}(2\pi; \lambda)$, too, has exactly $N$ zeros in the above half-plane).

Furthermore, for each $n > N$, $v(2\pi; \lambda)$ has exactly one simple zero in the egg-shaped region
\begin{equation}
\left|\sqrt{\lambda} - \frac{n}{2}\right| < \frac{1}{4}
\label{T28e}
\end{equation}
and $v(2\pi; \lambda)$ has no other zeros in the above region \cite{P-T}.

\section{A trace-like formula}
Let $\mu_1, \mu_2, \ldots$ be the zeros of $v(2\pi; \lambda)$ (counting multiplicities) labeled so that $|\mu_1| \leq |\mu_2| \leq \cdots$. Then, assuming that
$q \in C^2$ with
\begin{equation}
\langle q\rangle = \frac{1}{2\pi}\int_0^{2\pi} q(\xi)d\xi = 0,
\label{A0}
\end{equation}
we have the formula
\begin{equation}
\lim_n \sum_{j \leq n} \left( \mu_j - \frac{j^2}{4}\right) = \sum_{n = 1}^{\infty} \left( \mu_n - \frac{n^2}{4}\right) = -\frac{q(0)}{2}.
\label{A1}
\end{equation}
In the case of a real potential $q(x)$, where the zeros of $v(2\pi; \lambda)$ are simple and coincide with the Dirichlet eigenvalues of $H$ in the interval
$(0, 2\pi)$, such formulas are well known (see, e.g., the classical reference \cite{G-L}, which, however, contains a minor misprint regarding the sign in the
trace formula).

All proofs of trace formulas like \eqref{A1} that we have seen make use of the self-adjointness and, hence, are valid only for a real-valued $q(x)$. For this reason,
we have included the proof below, which works for complex potentials as well.

\medskip
 
\textit{Proof of formula} \eqref{A1}. The proof is done by estimating the contour integrals
\begin{equation}
\frac{1}{2\pi i} \oint_{C_n} \lambda
\left[\frac{\partial_{\lambda} v(2\pi; \lambda)}{v(2\pi; \lambda)} - \frac{\partial_{\lambda} \tilde{v}(2\pi; \lambda)}{\tilde{v}(2\pi; \lambda)} \right] d\lambda,
\label{A2}
\end{equation}
where $C_n$, $n \geq 1$ is the circle of radius $\left(\frac{n}{2} + \frac{1}{4}\right)^2$, centered at $0$, while $\partial_{\lambda}$ denotes 
the derivative with respect to $\lambda$.

Notice that, for $n$ sufficiently large the integral in \eqref{A2} is equal to the sum
\begin{equation*}
\sum_{j \leq n} \left( \mu_j - \frac{j^2}{4}\right).
\end{equation*}

To estimate the integrand of the contour integrals of \eqref{A2}, we begin with the asymptotic formula \eqref{T28c}. Dividing by $\tilde{v}(2\pi; \lambda)$
(recall \eqref{T1}) yields
\begin{equation}
m(\lambda) := \frac{v(2\pi; \lambda)}{\tilde{v}(2\pi; \lambda)} =
1 + \frac{q(0)}{2 \lambda}
+ O\left(\frac{1}{\lambda^{5/2}}\right),
\qquad
\lambda \to \infty,
\quad
\lambda \in \bigcup_{n=1}^{\infty} T_n,
\label{A3}
\end{equation}
where $T_n$, $n = 1, 2, \ldots$, are the annuli 
\begin{equation*}
T_n = \left\{\lambda \in \mathbb{C} : \left|\lambda - \left(\frac{n}{2} + \frac{1}{4}\right)^2 \right|  < 1 + n^{\alpha} \right\}
\end{equation*}
for some fixed $\alpha \in (0, 1)$. We restrict $\lambda$ in the annuli $T_n$, $n = 1, 2, \ldots$, so that the denominator $\tilde{v}(2\pi; \lambda)$ stays safely
away from $0$.

Notice that the asymptotic formula \eqref{A3} also implies
\begin{equation}
\frac{\tilde{v}(2\pi; \lambda)}{v(2\pi; \lambda)} =
1 - \frac{q(0)}{2 \lambda}
+ O\left(\frac{1}{\lambda^{5/2}}\right),
\qquad
\lambda \to \infty,
\quad
\lambda \in \bigcup_{n=1}^{\infty} T_n.
\label{A4}
\end{equation}
Next, let  $\Gamma \subset T_n$ be the circle of radius $n^{\alpha}$, centered at an arbitrary but fixed $\lambda \in C_n$. Then, Cauchy's 
integral formula together with \eqref{A3} and \eqref{A4} give
\begin{align}
m'(\lambda) 
&= \frac{\partial_{\lambda}v(2\pi; \lambda) \tilde{v}(2\pi; \lambda) - v(2\pi; \lambda) \partial_{\lambda} \tilde{v}(2\pi; \lambda)}{\tilde{v}(2\pi; \lambda)^2}
\nonumber
\\
&= \frac{1}{2\pi i} \oint_{\Gamma} \frac{m(z)}{(z - \lambda)^2}\, dz
= - \frac{q(0)}{2 \lambda^2} + o\left(\frac{1}{\lambda^{5/2}}\right),
\qquad
\lambda \to \infty,
\
\lambda \in \bigcup_{n=1}^{\infty} C_n.
\label{A5}
\end{align}
Finally, since
\begin{equation*}
\frac{\partial_{\lambda} v(2\pi; \lambda)}{v(2\pi; \lambda)} - \frac{\partial_{\lambda} \tilde{v}(2\pi; \lambda)}{\tilde{v}(2\pi; \lambda)}
= \frac{\tilde{v}(2\pi; \lambda)}{v(2\pi; \lambda)} \cdot\frac{\partial_{\lambda}v(2\pi; \lambda) \tilde{v}(2\pi; \lambda) - v(2\pi; \lambda) \partial_{\lambda} \tilde{v}(2\pi; \lambda)}{\tilde{v}(2\pi; \lambda)^2}
\end{equation*}
we get from the asymptotic formulas \eqref{A4} and \eqref{A5} that
\begin{equation*}
\lambda \left[\frac{\partial_{\lambda} v(2\pi; \lambda)}{v(2\pi; \lambda)} - \frac{\partial_{\lambda} \tilde{v}(2\pi; \lambda)}{\tilde{v}(2\pi; \lambda)}\right]
= - \frac{q(0)}{2 \lambda} \, + \, o\left(\frac{1}{\lambda^{3/2}}\right),
\quad
\lambda \to \infty,
\
\lambda \in \bigcup_{n=1}^{\infty} C_n.
\end{equation*}
Therefore,
\begin{equation*}
\sum_{j \leq n} \left( \mu_j - \frac{j^2}{4}\right) = \frac{1}{2\pi i} \oint_{C_n} \lambda
\left[\frac{\partial_{\lambda} v(2\pi; \lambda)}{v(2\pi; \lambda)} - \frac{\partial_{\lambda} \tilde{v}(2\pi; \lambda)}{\tilde{v}(2\pi; \lambda)} \right] d\lambda
= - \frac{q(0)}{2} + o\left(1\right)
\end{equation*}
as $n \to \infty$.
\hfill $\blacksquare$

\section{The system of equations for the $\mu$'s}
Suppose $q(x)$ is a real $C^3$ potential and $\mu_1(0), \mu_2(0), \ldots$ are the zeros of $v(2\pi; \lambda)$ associated with $q(x)$. Then \cite{Tr} the system of 
equations
\begin{equation}
\frac{d\mu_n}{d\xi} = \frac{n^2 \sqrt{\Delta(\mu_n)^2 - 4}}{8\pi \prod_{j \ne n} \left(\frac{\mu_j - \mu_n}{j^2/4}\right)},
\qquad
n = 1,2, \ldots,
\label{A6}
\end{equation}
where $\Delta(\lambda)$ is the Hill discriminant associated with $q(x)$, has a unique solution $\mu_1(\xi), \mu_2(\xi), \ldots$. Furthermore, under the appropriate 
choice of the signs of the square roots $\sqrt{\Delta(\mu_n)^2 - 4}$, the solution $\mu_1(\xi), \mu_2(\xi), \ldots$ of the system \eqref{A6} is the
set of zeros of $v_{\xi}(2\pi; \lambda)$, where $v_{\xi}(x; \lambda)$ is the solution of $H_{\xi} y = \lambda y$, where $H_{\xi}$ is the Hill operator associated with
$q_{\xi}(x) = q(x + \xi)$, satisfying $v_{\xi}(0; \lambda) = 0$ and $v_{\xi}'(0; \lambda) = 1$.

The derivation of the system of equations \eqref{A6} presented in \cite{Tr} remains valid for the case of a complex $q(x) \in C^3$.

\section{Meromorphic potentials}
As we have seen, if $\sigma(H) = [0, \infty)$, then $\Delta(\lambda) = 2\cos\left(2\pi \sqrt{\lambda}\right)$, hence
$\Delta(\lambda)^2 - 4 = -4\sin^2\left(2\pi \sqrt{\lambda}\right)$. Therefore, the system of equations \eqref{A6} takes the form
\begin{equation}
\frac{d\mu_n}{d\xi} = \sigma_n \frac{i n^2 \sin\left(2\pi \sqrt{\mu_n}\right)}{4\pi \prod_{j \ne n} \left(\frac{\mu_j - \mu_n}{j^2/4}\right)},
\qquad
n = 1,2, \ldots,
\label{A7}
\end{equation}
where $\sigma_n = \pm 1$.

Suppose now that we have coexistence for all $\lambda \neq 0, m^2/4$, where $m > 0$ is a given integer, while for $\lambda = m^2/4$ we do not have coexistence.
Then, $v_{\xi}(2\pi; n^2/4) = 0$  for all $n \geq 1$, $n \neq m$. Consequently, $\mu_n(\xi) = n^2/4$ for all  $n \geq 1$, $n \neq m$ and the system \eqref{A7}
reduces to a single differential equation for $\mu_m(\xi)$:
\begin{equation*}
\frac{d\mu_m}{d\xi} = \pm 2 i \sqrt{\mu_m} \left(\dfrac{m^2}{4} - \mu_m\right).
\end{equation*}

This equation can be easily solved and
from its solutions we can obtain the associated potentials $q(x)$ (via formula \eqref{A1}), which turn out to be the meromorphic Gasymov potentials (recall our Example)
\begin{equation*}
q_m(x) = \frac{2 m^2 a e^{imx}}{(a e^{imx} + 1)^2} = \frac{2 m^2 a^{-1} e^{-imx}}{(a^{-1} e^{-imx} + 1)^2},
\qquad
a \neq 0, \quad |a| \neq 1.
\end{equation*}
Applying successive Darboux transformations we can obtain meromorphic potentials which are
not Gasymov but whose spectrum is $[0, \infty)$. This fact was first suggested in R. Carlson's
paper \cite{Car}. Hence, our conjecture is not true for meromorphic potentials. Notice that each Darboux transformation destroys the coexistence of one periodic
or antiperiodic eigenvalue.

\section{Another analog of Borg's theorem}
\textbf{Theorem.} Suppose $q \in C^2$ and $\sigma(H) = [\,0, \infty)$, thus $\Delta(\lambda) = 2\cos(2\pi \sqrt{\lambda})$. Furthermore, suppose that we have coexistence at
$\lambda = n^2/4$, for every integer $n \geq 1$. Then $q(x) \equiv 0$.

\smallskip

\textit{Proof}. Notice that coexistence at $\lambda = n^2/4$, for every integer $n \geq 1$, implies that both $u(x; n^2/4)$ and  $v(x; n^2/4)$ are
Floquet solutions and, consequently, periodic or antiperiodic, since $\rho(n^2/4) = \pm1$. Therefore, $v(2\pi; n^2/4) = 0$ for every integer $n \geq 1$.
From the asymptotic formulas \eqref{T28d} and \eqref{T28e} it follows that 
these are the only zeros of $v(2\pi; \lambda)$ and that all these zeros are simple. Thus, the zeros of $v(2\pi; \lambda)$ (counting multiplicities) are
\begin{equation*}
\mu_n = \frac{n^2}{4},
\qquad
n \geq 1.
\end{equation*}
Furthermore, the same is true for the quantity $v_{\xi}(2\pi; \lambda)$ associated with the shifted operator $H_{\xi}$, for every $\xi \in \mathbb{R}$.

Therefore, in view of formula \eqref{A1}, we get
\begin{equation*}
0 = \sum_{n = 1}^{\infty} \left( \frac{n^2}{4} - \frac{n^2}{4}\right)
= \sum_{n = 1}^{\infty} \left[ \mu_n(\xi) - \frac{n^2}{4}\right] = -\frac{q_{\xi}(0)}{2} \equiv -\frac{q(\xi)}{2},
\end{equation*}
i.e. $q(\xi) \equiv 0$.
\hfill $\blacksquare$

\medskip

\textbf{Acknowledgment.} The author wants to thank Professors Tuncay Aktosun, Ricardo Weder, and the other organizers of the
\textit{Analysis and Mathematical Physics 2024} (online) Conference for inviting him to participate with a talk. The present article is an adaptation of the author's presentation.

\end{document}